\documentclass[12pt]{article}
\usepackage{amsmath,amsfonts,amsthm}

\def\hybrid{\topmargin 0pt      \oddsidemargin 0pt
        \headheight 0pt \headsep 0pt
        \voffset=-0.5cm
        \textwidth 6.5in        
        \textheight 9in         
        \marginparwidth 0.0in
        \parskip 5pt plus 1pt   \jot = 1.5ex}

\hybrid

\def\beq{\begin{equation}}
\def\eeq{\end{equation}}
\def\bea{\begin{eqnarray}}
\def\eea{\end{eqnarray}}

\def\G{\Gamma}
\def\g{\gamma}
\def\s{\sigma}

\def\C{{\cal C}}

\def\a{\alpha}

\def\R{{\cal R}}
\def\A{{\cal A}}

\def\D{{\cal D}}
\def\F{{\cal F}}

\def\L{{\cal L}}

\def\E{{\cal E}}

\def\P{{\cal P}}

\def\SP{{\cal S}}

\def\res{{\rm res}}
\def\T{{\rm T}}

\def \matrix #1 {\left(\begin{array}{cc} #1 \end{array}\right)}

\newtheorem{theo}{Theorem}[section]

\newtheorem{cor}[theo]{Corollary}
\newtheorem{lem}[theo]{Lemma}

\theoremstyle{definition}

\newtheorem{rem}[theo]{Remark}

\def\bbZ{{\mathbb Z}}
\def\bbC{{\mathbb C}}

\begin{document}

\title{Commuting difference operators and the combinatorial Gale transform}
\author{I.Krichever
\thanks{Department of Mathematics, Columbia University; High School of Economics, Kharkevich Institute for Information Transmission Problems, and Landau Institute for Theoretical Physics, Moscow, Russia.
Research is supported in part by Simons foundation and in part by RFFI grants  13-01-12469 and 14-01-00012-a}}

\date{}

\maketitle

\begin{abstract} We study the spectral theory of $n$-periodic strictly triangular difference operators $L=T^{-k-1}+\sum_{j=1}^k a_i^j T^{-j}$ and the spectral theory of the "superperiodic" operators for which all solutions of the equation $(L+1)\psi=0$ are (anti)periodic. We show that for a superperiodic operator $L$ there exists a unique superperiodic operator $\L$ of order $(n-k-1)$ which commutes with $L$ and show that the duality $L\leftrightarrow \L$ coincides up to a certain involution with the combinatorial Gale transform recently introduced in \cite{tabachnikov1}
\end{abstract}

\section{Introduction}

Recently, the theory of linear difference equations of the form
\beq\label{oper}
V_i=a_i^1V_{i-1}-a_i^2V_{i-2}+\cdots+(-1)^{k-1} a_i^{\,k}V_{i-k}+(-1)^{k}V_{n-k-1}
\eeq
with real periodic coefficients
\beq\label{coef}
a_i^j=a_{i+n}^j\in \mathbb R\,,\ i\in \mathbb Z,\ 1\leq j\leq k\,,
\eeq
whose {\it all} solutions are $n$-(anti)periodic
\beq\label{vi}
V_{i+n}=(-1)^k V_i
\eeq
has seen a burst of interest due to connections with the theory of $n$-gons in projective space  (extensively studied in \cite{gale,gelfand,kapranov,eisenbud}) and $SL_{k+1}$-frieze patterns (created by Coxeter and Conway in \cite{coxeter1,coxeter2} and developed later in \cite{propp,assem,bergeron,morier1,morier2,ovsienko1,sophie}). The latter is a particular case of more general patterns arising in the theory of {\it bilinear discrete Hirota} equations (see \cite{zabrodin}, and references therein).

Our study was motivated by recent work \cite{tabachnikov1}, where it was shown that the moduli space $\E_{k+1,\,n}$ of equations (\ref{oper}) satisfying constraint (\ref{vi}) is isomorphic, as algebraic variety, to the moduli space $\F_{k+1,\,n}$ of $SL_{k+1}$-frieze patterns of width $w=n-k-2$, and that, if $n$ and $k+1$ are coprime, then both spaces are isomorphic to the moduli space $\C_{k+1,n}$ of $n$-gons in $\mathbb{RP}^k$. The main result of \cite{tabachnikov1} is a description of the {\it combinatorial Gale} transform, which is the duality between the spaces $\E_{k+1,\,n}$ and $\E_{n-k-1,\,n}$ analogous to the classical Gale duality on the moduli spaces of point configurations. The former and the latter are induced by the duality of the Grassmanianns $Gr(k+1,n)$ and $Gr(n-k-1,n)$.

Recall, that if $\widehat V$ is  $(n\times (k+1))$ matrix of rank $k+1$ representing a point of
$Gr(k+1,n)$, then the dual point is represented by $(n\times (n-k-1))$ matrix $\widehat W$ such that $\widehat W^T\,\widehat V$=0. The space of $n$-gons in $\mathbb P^k$ is a quotient of the Grassmaniann by the torus action,
$\C_{k+1,n}\simeq Gr(k+1,n)/\mathbb T^{n-1}$ (\cite{gelfand}). The points of the Grassmanianns representing dual $n$-gons
are defined by the equation $W^TDV$, where $D$ is a non-degenerate diagonal $(n\times n)$ matrix.

If $V^{(j)}$ is a basis of solution of equation (\ref{oper}) then the matrix $\widehat V_{(i)}:=\left(V_{i+1}^{(j)},V_{i+2}^{(j)},\ldots, V_{i+n}^{(j)}\right)$ defines an imbedding of $\E_{k+1,n}$
into the Grassmaniann $Gr(k+1,n)$ (depending  on the choice of $i$). As shown in \cite{tabachnikov1} (Proposition 4.3.1) the defining property of the combinatorial Gale dual equations is the relation:
\beq\label{matrixdual}
\widehat W_{(i)}^{\,T}\ \widehat \varepsilon\  \widehat V_{(j)},\ \ i-j=n-k-1\  ({\rm mod}\  n)\,,
\eeq
where $\widehat \varepsilon$ is the diagonal matrix $(1,-1,1,\ldots)$.

The simplest example of the combinatorial Gale transform is related to Gauss' {\it pentagramma  mirificum}
\cite{gauss}. It is a duality between $5$-periodic difference equations of third and second order. More precisely, the duality states that if all solutions of the equation
\beq\label{mir3}
V_i=a_i V_{i-1}-b_i V_{i-2}+V_{i-3},\ a_i=a_{i+5},\ b_i=b_{i+5},
\eeq
are $5$-periodic, then all solutions of the equation
\beq\label{mir2}
W_i=a_i W_{i-1}-W_{i-2}
\eeq
are $5$-antiperiodic, $W_i=-W_{i+5}$. In terms of $n$-gons this transformation sends projectively equivalence classes of pentagons in $\mathbb P^2$ to those in $\mathbb P^1$. In terms of frieze patterns it is a duality between 5-periodic Coxeter friezes and 5-periodic $SL_3$-friezes.

The first goal of this work is to introduce a {\it spectral parameter} in equation (\ref{oper}) and to develop
the Block-Flouque spectral theory of the corresponding difference operators. We would like to stress, that the idea to introduce a spectral parameter in (\ref{oper}) is not new by itself. In \cite{tabachnikov} the spectral parameter $s$ was introduced in equation (\ref{mir3}) as follows:
\beq\label{mir3s}
V_i=sa_i V_{i-1}-s^{-1}b_i V_{i-2}+V_{i-3},\ a_i=a_{i+5},\ b_i=b_{i+5},
\eeq
and used for the construction of a complete set of integrals in involution for the pentagram map defined in
\cite{schwartz}. Algebraic-geometrical integrability of the pentagram map was proved in \cite{soloviev}
using the Lax operator with the spectral parameter which is gauge equivalent to (\ref{mir3s})
(for Liouville-Arnold integrability see \cite{ovsienko}). In a similar way
the spectral parameter was introduced in higher order linear equations for the proof of algebraic geometrical integrability of higher pentagram maps (\cite{khesin-sol,khesin-sol1}) (see also \cite{maria}). An  approach to introduce the spectral parameter in equation (\ref{oper}) using cluster algebras was proposed in \cite{fock}, \cite{gekhtman}, which for $SL_3$-case turned out to be equivalent (\ref{mir3s}).

The spectral parameter $s$ in (\ref{mir3s}) reflects the scaling invariance of the pentagram map (see \cite{tabachnikov}). In this paper we consider the {\it different} spectral problem. Namely, let $\D_{k+1,\,n}$ be the (affine) space of monic "strictly" triangular $n$-periodic difference operators with complex (or real) coefficients:
\beq\label{Lax}
\left\{L\in \D_{k+1,n}\ |\ L=T^{-k-1}+\sum_{j=1}^k a_i^j T^{-j}\,,\ \  a_i^j=a_{i+n}^j\right\}
\eeq
Here and below, $T$ is the shift operator acting on an infinite sequence $\psi=\{\psi_i\},\, i\in \mathbb Z$, i.e., $(T\psi)_i=\psi_{i+1}$. It will be also assumed, unless stated otherwise, that the leading coefficient of $L$
does not vanish,
\beq\label{a1}
a_i^1\neq 0
\eeq
The equation for eigenvectors of $L$
\beq\label{eigenL}
L\psi=E\psi
\eeq
after the gauge transformation
\beq\label{gauge}
V=\varepsilon\psi\ \ \ \varepsilon:=\{\varepsilon_i=(-1)^{i}\}\,,
\eeq
coincides with (\ref{oper}) for $E=-1$.
\begin{rem} Similarly, the spectral parameter in equation(\ref{oper}) can be introduced as the coefficient
at the last term of the right hand side of the equation, but two ways of introducing the spectral parameter are in fact equivalent due to the involution on $\D_{k+1,n}$:
\beq\label{inv}
L\longmapsto L^\s:=T^{-k-1}(L^*+1)-1,
\eeq
where $L^*$ is the formal adjoint operator
\beq\label{adj}
L^*=T^{k+1}+\sum_{j=1}^k T^ja_i^j=T^{k+1}+\sum_{j=1}^k a_{i+j}^jT^j.
\eeq
\end{rem}
\medskip
Under the gauge transformation the constraint (\ref{vi}) defining the special triangular  operators $L\in \E_{k+1,\, n}\subset \D_{k+1,n}$ takes the form:
{\it all solutions $\Psi$ of the equation $(L+1)\Psi=0$ have the following monodromy property}:
\beq\label{psiper}
\Psi_i=(-1)^{n+k}\Psi_{i-n}
\eeq
For brevity, throughout the rest of the paper the operators $L\in \E_{k+1,n}$ will be called {\it superperiodic}.

\begin{rem} Note, that the choice of the multiplicator $\mu=(-1)^{n+k}$ in (\ref{psiper}) can be seen as a choice of the normalization in the following slightly more general setting.
For a pair of complex numbers $(e,\, \mu)$ consider the (affine) space  of triangular $n$-periodic operators of order $(k+1)$  such that {\it all solutions} of the equation $(L-e)\Psi=0$ are the Bloch solutions with the same multiplicator $\mu$,
\beq\label{special}
\E_{k+1,\,n}^{e,\,\mu}=\left\{L\in \D_{k+1,n} |\, (L-e)\Psi=0 \Rightarrow \mu\Psi_i=\Psi_{i-n}\right\}
\eeq
The same arguments as in \cite{tabachnikov1} show that this space is nonempty if and only if
\beq\label{ab}
\mu^{k+1}=(-1)^{kn}e^n\,.
\eeq
Indeed, the operator $L$ is monic. Therefore, if $\Psi^{(j)}$ is the basis in the space of solutions of (\ref{special}), then $\det \Phi_{(i)}=(-1)^k e^{-1} \det \Phi_{(i-1)}$, where $\Phi_{(i)}:=\left\{\Psi^{(j)}_{i-s}\right\}^{j=0,\ldots, k}_{s=0,\ldots,k}$. The monodromy operator $T^{-n}$ restricted onto the space spanned by $\Psi^{(j)}$ is the scalar $\mu$. Therefore, the latter equation implies (\ref{ab}).

Over the complex numbers it is obvious that all admissible pairs $(e,\mu)$ are in the same orbit of the scaling transformation $(e,\mu)\to (c^{k+1}e,c^n\mu),\ c\in \mathbb C$. Over the real it can be directly verified that if $n$ and $(k+1)$ are co-prime, then the parity of $nk$ coincides with the parity of $k$ and {\it any real admissible pair can be obtained by the real scaling from the pair} $(-1,(-1)^{n+k})$ corresponding to the case of superperiodic operators.

The scaling transformation above reflects the scaling transformation of operators
\beq\label{scaling}
\tau_c: \D_{k+1,n}\longmapsto \D_{k+1,n},\ \tau_c(L)= c^{k+1} CLC^{-1}\,,
\eeq
where $C$ is the diagonal operator, $(C\psi)_i=c^i\psi_i$. Notice, that the multiplication by $c^{k+1}$ of the gauge transformed operator $CLC^{-1}$ is introduced in order to make the operator $\tau_c(L)$ monic.
\end{rem}

In the theory of ordinary linear differential operators it is well-known that if the kernel of an operator $B$ is a subspace of the kernel of another operator $A$, $ker B\subset ker A$, then the latter is divisible by the former from the right and from the left, i.e. there exist operators $C_l,C_r$ such that $A=BC_r=C_lB$. Analogous to the the case of differential operators, it is easy to show, that the original defining characterization of the superperiodic operators is equivalent to the following one:

\begin{lem}\label{lemdivis} Let $L\in \D_{k+1,n}$ be a $n$-periodic triangular operator of order $k+1$. Then it is
superperiodic if and only if the operator $T^{-n}-(-1)^{n+k}$ is divisible from the right and from the left by $L+1$, i.e. there exist operators $L_r,L_l$ in $\D_{n-k-1,\,n}$ such that the equations
\beq\label{divis}
T^{-n}-(-1)^{k+n}=(L_l-(-1)^{k+n})(L+1)=(L+1)(L_r-(-1)^{k+n})
\eeq
hold.
\end{lem}

In section 3 we show that if the period $n$ and the order $k+1$ of $L$ are co-prime, then the Bloch solutions of (\ref{eigenL}), i.e. solutions that are eigenvectors for the monodromy operator $T^n$:
\beq\label{monod}
w\psi_i=\psi_{i-n},
\eeq
coincide with the discrete analog of the Baker-Akhiezer function defined on the spectral curve of $L$.
This Bloch spectral curve is different from the spectral curves which have arisen in the theory of the pentagram map but is relevant to the combinatorial Gale transform.

More precisely, using the identification of the Bloch solutions of (\ref{Lax}) with the discrete Baker-Akhiezer function we establish a connection of the spectral theory of triangular operators with the theory of rank 1 commuting difference operators developed in \cite{kr-diff,mumford}, and through this connection we clarify the defining  constraint of operators $L\in \E_{k+1,n}$ and the combinatorial Gale duality.

An appearance of commuting operators in the theory of superperiodic difference operators is the direct corollary of our first main result:
\begin{theo}\label{main} If $n$ and $k+1$ coprime, then:

(i) the operators $L_r$ and $L_l$ defined by equation (\ref{divis}) for a superperiodic operator $L$ coincide, i.e. $L_r=L_l$;

(ii) under involution (\ref{inv}) the operator $\L:=L_r=L_l$ gets transformed to the operator which is gauge equivalent to the operator ${\cal G}(L)$ Gale dual to $L$:
\beq\label{duality}
{\cal G}(L)=\varepsilon\, \L^\s \varepsilon^{-1}
\eeq
\end{theo}
The statement $(i)$ of the theorem and equations (\ref{divis}) imply:
\begin{cor} The operators $L$ and $\L$ commute with each other:
\beq\label{com25}
[L,\L]=0.
\eeq
\end{cor}
\noindent
As we see, the existence of a commuting operator is the necessary condition for a periodic triangular operator to be superperiodic. It turns out that this condition is necessary and sufficient for operators that can be obtained from the superperiodic operators by the scaling transformation (\ref{scaling}).

\begin{theo}\label{main1} Suppose that for an operator  $L\in \D_{k+1,n}$ with $2(k+1)<n,\ (n,k+1)=1,$ there exists a commuting operator $K\in \D_{n-k-1,\,n}$. Then in general position:

(i) the operator $L$ is superperiodic up to the scaling transformation (\ref{scaling}) for some constant $c$, i.e., $\tau_c(L)\in \E_{k+1,n}$;

(ii) there is a unique polynomial $P$ such that the operator $\L=K+P(L)$
is superperiodic up to the scaling transformation (\ref{scaling}).
\end{theo}
\noindent
Here the "general position" means: for a Zariski open set of parameters. The precise meaning will be given in the proof of the theorem in Section 4.

\section{The discrete Baker-Akhiezer function}

To begin with let us recall basic facts of the theory of commuting rang 1 difference operators.

Let $\G$ be a smooth algebraic curve of genus $g$. Fix two points
$p_{\, \pm}\in \G$, and let $D=\g_1+\cdots+\g_{g}$ be a generic effective divisor on $\G$ of
degree $g$.  By the Riemann-Roch theorem one computes $h^0(D+i(p_+-p_-))=1$, for any $i\in\bbZ$, and for
$D$ generic. We denote by $\psi_{i}(p),\ p\in \G,$ the unique section of this bundle. This means that $\psi_{i}$ is the unique up to a constant factor meromorphic function such that (away from the marked points $p_{\pm}$) it has poles only at $\g_s$, of multiplicity not greater than the multiplicity of $\g_s$ in $D$, while at the point $p_+$ (resp. $p_-$) the function $\psi_{i}$ has zero (resp. pole) of order $i$.

If we fix local coordinates $z$ in the neighborhoods of marked points, then the Laurent series for $\psi_{i}(p)$, for $p\in\G$ near a marked point, has the form
\beq\label{1}
\psi_{i}=z^{\pm i}\left(\sum_{s=0}^{\infty}\xi_s^{\,
\pm}(i)z^{s}\right), \ \ z=z(p), \
p\to p_{\, \pm}
\eeq
The function $\psi:=\psi_i(p)$ of the discrete variable $i$ and the point $p$ of the curve $\G$ is called the discrete Baker-Akhiezer function (\cite{kr-diff}). Throughout the paper it will be assumed that $\psi$ is normalized by the condition $\xi^+_0=1$. Notice that under the change of local coordinate $z\to \tilde z=cz+O(z^2)$ near the marked point $p_+$ the normalized Baker-Akhiezer function gets transformed to
\beq\label{trans}
\tilde \psi_i(p)=c^i \psi_i(p)
\eeq
(compare with the scaling transformation (\ref{scaling}).

Let $\A(p_+,p_-)$ be the ring of meromorphic functions on $\G$ that are holomorphic away from the marked points $p_{\pm}$. From the uniqueness of the Baker-Akhiezer function it easy follows that
\begin{lem}[\cite{kr-diff}] \label{comring} For each  $A\in \A(p_+,p_-)$ there exists a unique difference operator
\beq\label{La}
L_A=\sum_{j=-k_-}^{k_+} a_i^j T^{-j}  \ \ \ \Leftrightarrow \ \ \ (L_A\psi)_i=\sum_{j=k_-}^{k_+} a_i^j \psi_{i-j}
\eeq
such that
\beq\label{Lax11}
L_A\psi=A\psi\,.
\eeq
Here $k_{\pm}$ are the orders of poles of $A$ at $p_{\pm}$.
\end{lem}
\begin{cor} The operators $L_A$ commute with each other, i.e. $[L_A,L_B]=0$.
\end{cor}
The coefficients of the operator $L_A$ are defined recurrently from the system of equations obtained by the substitution of the expansion of $\psi$ and $A$ near the marked points. For example, if $A=c_0z^{-k_+}+c_1z^{-k_+1}+\cdots$ is the expansion of $A$ near the points $p_{+}$, then
\beq\label{coefficients}
a_i^{k_+}=c_0,\ \  a_i^{k_+-1}=c_0(\xi_1^+(i+k_+)-\xi_1^+(i))-c_1\,\ldots
\eeq
Similarly, if $A=c_0^- z^{-k_-}+\cdots$ is the expansion of $A$ near $p_-$, then
\beq\label{coefficients-}
a_i^{-k_-}=\frac{c_0^-\xi_0^-(i)}{\xi_0^-(i+k_-)},\ \ldots
\eeq
\begin{rem}
For further use note, that if $A$ is in the subring $\A_+(p_+,p_-)\subset\A(p_+,p_-)$ of meromorphic functions on $\G$ that have pole only at $p_+$ and vanish at $p_-$, then the operator $L_A$ is strictly lower triangular.
In that case the formula (\ref{coefficients-}) for the leading coefficient $a_i^1$ takes the form
\beq\label{coeff-a1}
a_i^1=e^{x_i-x_{i-1}+l}, \ \ x_n:=\ln \xi_0^-(i), \ \ l:=\ln c_0^-.
\eeq
\end{rem}

\begin{rem} The correspondence above extends to the case of singular curves. More precisely, (\cite{kr-diff,mumford}):
there is a natural correspondence
\beq\label{corr}
\A\longleftrightarrow \{\G,p_{\pm},  \F\}
\eeq
between commutative rings $\A$ of ordinary linear
difference operators containing a pair of monic operators of co-prime orders, and
sets of algebraic-geometrical data $\{\G,p_{\pm}, [z]_1, \F\}$, where $\G$ is an
algebraic curve with a fixed first jet $[z]_1$ of a local coordinate $z$ in the neighborhood of a smooth
point $P_+\in\G$ and $\F$ is a torsion-free rank 1 sheaf on $\G$ such that
\beq\label{sheaf}
h^0(\G,\F(ip_+-ip_-))=h^1(\G,\F(ip_+-ip_-))=0.
\eeq
The correspondence becomes one-to-one if the rings $\A$ are considered modulo conjugation
$\A'=g\A g^{-1}, g=\{g_i\}$ (in its final form correspondence (\ref{corr}) is due to Mumford) .
\end{rem}
The discrete Baker-Akhiezer function, and therefore, the coefficients of the corresponding difference operators can be expressed in terms of the Riemann theta-function
\beq\label{theta}
\theta(z):=\theta(z|B)=\sum_{m\in \mathbb Z^g}e^{2\pi i (m,z)+\pi i (Bm,m)}\,,\ \ z\in \mathbb C^g\,,
\eeq
defined by the matrix $B$ of $b$-periods of the normalized holomorphic differentials $d\omega_j$ on $\G$
\beq\label{matrixB}
B_{j,k}=\oint_{b_k}d\omega_j\,,\ \ \ \delta_{j,k}=\oint_{a_k}d\omega_j
\eeq

\begin{lem}[\cite{kr-diff}] The Baker-Akhiezer function is given by the formula
\beq
\psi_i(p)=
\frac{\theta(A(p)+i U+Z) \theta (A(p_+)+Z)}
{\theta(A(p_+)+i U+Z)\theta(A(p)+Z)} \ e^{i\Omega(p)}, \label{psitheta}
\eeq
Here
a) $\Omega_0(p)$ is the abelian integral,
$
\Omega_0(p)=\int^p d\Omega_0,
$
corresponding to the unique normalized,
$
\oint_{a_k} d\Omega_0=0,
$
meromorphic differential on $\Gamma$, which has simple poles at the marked point $p_{\pm}$ with residues
$\pm 1$, respectively;

b) $A(p)$ is the Abel transform, i.e., a vector with the coordinates
$
A(p)=\int^p d\omega_k\,;
$

c) $U$ is the vector $A(p_+)-A_(p_-)$

d) $Z$ is an arbitrary vector (it corresponds to the divisor of poles of
Baker-Akhiezer function).
\end{lem}
From (\ref{psitheta}) it easy follows that the coefficients of the commuting difference operators $L_A$ are
{\it quasi-periodic} functions of the variable $i$. The operators with periodic coefficients are singled  out by the following constraint.
\begin{lem}\label{periodic} Suppose that on $\G$ there exists a meromorphic function $w=w(p)$ which is
holomorphic away of the marked point $p_+$ where it has a pole of order $n$, and which has zero of order $n$ at the marked point $p_-$. Then the coefficients of the operators $L_A$ are $n$-periodic.
\end{lem}
\noindent
{\it Proof.} The function $w$, if exists, is unique up to multiplication by a constant. If the first jet of a local coordinate near the point $p_+$ is fixed, then $w$ can be normalized by the condition $w=z^{-n}+O(z^{-n+1})$.
The uniqueness of $\psi_i$ implies then, that equation (\ref{monod}) holds. Hence, $L_A$ are $n$-periodic.

\bigskip

{\bf Real operators.} The formula (\ref{psitheta}) for a generic $i\in C$ defines $\psi_i(p)$ as a multivalued function on $\G$. A single-valued branch of it can be defined on $\G$ with a cut between the marked points which does not intersect the chosen basis of $a$- and $b$-cycles. The coefficients $a_i^j$ of the corresponding operators can be seen as {\it meromorphic} functions of the complex variable $i$.
The spectral data corresponding to operators with real (for $i\in \mathbb R$) coefficients are singled out as follows:

\begin{lem} Suppose that on an  algebraic curve $\G$ with two marked points $p_{\pm}$ there exists an antiholomorphic involution $\tau$ for which the marked points are fixed, $\tau(p_{\pm})=p_\pm$. Then the Baker-Akhiezer function $\psi$ corresponding to a real divisor $D=\tau(D)$ satisfies the relation
\beq\label{tau}
\psi(p)=\bar\psi(\tau(p))
\eeq
(provided that the coordinate $z$ in the neighborhood of $p_+$ used for the normalization of $\psi$ is also real, $\bar z=\tau^* z$).
\end{lem}
Equation (\ref{tau}) implies that if $A=\bar A(\tau(p))$ then the corresponding operators $L_A$ have real coefficients. Still they might be singular functions of the variable $i\in \mathbb R$.
In the framework of the finite-gap theory there are two basic types of conditions which are sufficient for regularity of the corresponding operators. Let us present one of them which is relevant to the case under consideration.

Recall, that an antiholomorphic involution $\tau$ of a genus $g$ smooth  algebraic curve $\G$ has at most $(g+1)$
fixed ovals. The curves having the maximal number of fixed ovals $A_0,\ldots,A_g$ are called {\it $M$-curves}.

\begin{lem}\label{mcurves} Let $\G$ be an $M$-curve with two marked points $p_{\pm}\in A_0$ and let $D$ be a set of $g$ points $\g_s\in A_s, s=1,\ldots,g$. Then for real functions $A=\bar A(\tau(p))\in \A(p_+,p_-)$ the coefficients of the operators $L_A$ corresponding to the Baker Akhiezer function defined by $D$ are real and nonsingular. Moreover,
their leading coefficients  $a_i^{-k_-}$   are sign definite:
\beq\label{signum}
a_i^{k_+}=c_0,\ \ \ {\rm sgn}\, (a_i^{-k_-})={\rm sgn} \,c_0^-
\eeq
\end{lem}
\noindent
The proof is standard in the finite gap theory (see for example the proof of the analogous statement for the difference Schr\"odinger operator in \cite{viniti}).

In what follows we will use one more basic concept of the algebraic-geometrical integration theory.

\medskip
\noindent{\bf The dual Baker-Akhizer function.}
The concept of the dual Baker-Akhiezer function $\psi^+$ is universal
and is at the heart of Hirota's bilinear form of soliton equations.
In the discrete case the dual Baker-Akhiezer function corresponds to the dual
divisor $D^+=\g_1^++\cdots+\g_g^+$, which is defined by the equation
\beq\label{D+D}
D+D^+=K+p_-+p_+\in J(\G)
\eeq
where $K$ is the canonical class. In other words: the points of $D$ and $D^+$ are zeros of the meromorphic differential $d\Omega$ with simple poles at $p_{\pm}$ with residues $\pm 1$, respectively. The
dual Baker-Akhiezer function is then defined by the following analytic properties:

i) the function $\psi_i^+$ (as a function of the variable $p\in \G$) is meromorphic everywhere except for the points $p_{\pm}$ and has at most simple poles at the points $\gamma_1^+,\ldots,\gamma_g^+$ (if all
of them are distinct);

(ii) in a neighborhood of the point $P_{\a}$ the function $\psi$ has the
form
\beq\label{dual1}
\psi_{i}^+=z^{\mp i}\left(\sum_{s=0}^{\infty}\chi_s^{\,
\pm}(i)z^{s}\right), \ \ z=z(p), \
p\to p_{\, \pm},\ \ \chi^+_0=1.
\eeq
In fact it is the same Baker-Akhiezer type function and, therefore, admits the same type of explicit theta-function formula:
\beq
\psi_i^+(p)=
\frac{\theta(A(p)-i U-Z-A(p_+)-A(p_-))\, \theta (A(p_-)+Z)}
{\theta(A(p_-)+i U+Z)\, \theta(A(p)-Z-A(p_+)-A(p_-))} \ e^{-i\Omega_0(p)}, \label{dualpsitheta}
\eeq
\begin{lem}\label{lem3.1} Let $\psi$ and $\psi^+$ be the Baker-Akhiezer function and its dual. Then the following
equations hold:
\beq\label{resbad}
\res_{p_+} (\psi_i^+\psi_j)\,d\Omega=\delta_{i,\,j}, \ \  i,j\in \mathbb Z.
\eeq
\end{lem}
By definition of the duality, the differential on the left-hand side of (\ref{resbad}) for $i\neq j$
has pole only at {\it one} of the marked points $p_{\pm}$. Hence, its residue vanishes. The differential $\psi_i^+\psi_i d\Omega$ has poles at $p_+$ and $p_-$ with residues $\pm 1$ respectively. Thus the lemma is proven.

\begin{cor}Let $\psi$ be the Baker-Akeiezer function and let $L_A$ be the linear operators
such that (\ref{Lax11}) hold. Then the dual Baker-Akhiezer function is a solution of the formal adjoint equation
\beq\label{2.9adj-2Toda}
\psi^+L_A=A\psi^+
\eeq
\end{cor}
Recall that the right action of a difference operator is defined as the formal adjoint action, i.e., $f^+T=T^{-1}f^+$.

\section{The Bloch-Floquet theory of triangular periodic difference operators}

Let us briefly recall the conventional setting of the spectral theory of periodic difference operators. Consider an $n$-periodic difference operator $L$. Then the monodromy operator $T^{-n}$ preserves
$(k+1)$-dimensional space $\L(E)$ of solutions of the equation $L\psi=E\psi$. Let $T^{-n}(E)$ be the restriction of $T^{-n}$ onto $\L(E)$. The common eigenvalues of $L$ and the monodromy operators satisfy the algebraic equation
\beq\label{curve}
R(w,E)=\det (w\cdot 1-T^{-n}(E))=0
\eeq
The same equation can be obtained if one introduces  $n$-dimensional space ${\cal T}(w)$ of solution of equation (\ref{monod}) and denote by $L(w)$ the restriction of $L$ onto $\T(w)$. Then
\beq\label{curve1}
R(w,E)=\det (L(w)-E\cdot 1)=0.
\eeq
Out first observation which turns out to be crucial for the further considerations is the following:
\begin{lem}\label{lem-curve} If $n$ and $k+1$ are coprime, then the spectral curve of an operator   $L\in \D_{k+1,\,n}$ is defined by the equation of the form:

\beq\label{R}
R(w,E)=w^{k+1}-E^n+\sum_{i>\,0,\, j\geq \,0,\ ni+(k+1)j<n(k+1)}  r_{ij} w^iE^j=0.
\eeq
with
\beq\label{r}
r:=r_{1,\, 0}=\prod_{j\,=1}^n a_j^1.
\eeq
\end{lem}
\begin{rem}\label{important}
The statement above is the statement that the Newton polygon of $R$
is a triangle with the vertices $(0,n),(k+1,0), (0,1)$ with the explicit form of the coefficients corresponding to the vertices. Before, presenting a proof of the lemma, let us elaborate on specifics of the  algebraic curves defined by equation (\ref{R}).

From (\ref{R}) it follows that the affine curve defined by equation $R(w,E)=0$
at infinity is compactified by one {\it smooth} point $p_+$\,. At this point the function $w=w(p),\ p\in \G,$ has pole of order $n$ and $E=E(p)$ has pole of order $k+1$. Notice, that (\ref{R}) implies also that if $w=0$, then $E=0$. Moreover, if $r\neq 0$, then one branch of the multivalued function $w(E)$ defined by (\ref{R}) near the origin is of the form
\beq\label{origin}
w=r^{-1}E^n\left(1+\sum_{s=1}^{\infty} v_sE^s\right)
\eeq
Hence, the curve $\G$ has another marked {\it smooth} point $p_-$ at which $w$ has zero of order $n$ and $E$ has simple zero. The degrees of zero and pole divisors are equal, therefore, $w$ does not vanish on $\G\setminus p_-$.

Let $\SP_{k+1,\,n}$ be the family of curves $\G$ defined by  equations of the form (\ref{R}). It is of dimension $\frac{k(n+1)}2$ (the number of the coefficients $r_{ij}$). For generic values of coefficients $r_{ij}$ the curve $\G$ is smooth and has genus $g=\frac{k(n-1)}2 $\,. Let $\P_{k+1,\,n}$ be the Jacobian bundle over an open subspace of
smooth curves in $\SP_{k+1,\,n}$. A generic pair $\{\G, D\}\in \P_{k+1,\,n}$ defines a unique Baker-Akhiezer function $\psi$, and then, by Lemma \ref{comring}, the operator $L_E$ corresponding to the function $E=E(p)$. This function has pole of order $k+1$ at $p_+$ and vanishes at $p_-$. Hence, $L_E$ is a strictly lower triangular. Moreover, corollary \ref{periodic} implies that $L_E$ is $n$-periodic. The established correspondence of an open domain
\beq\label{correspondence}
\P_{k+1,\,n}^0\subset \P_{k+1,\,n}\longmapsto \D_{k+1,\,n}
\eeq
will be referred as {\it the inverse spectral transform}.

\medskip
\noindent{\bf Example.} For $k=1, n=2m+1$ the family $\SP$ is the family of hyperelliptic curves defined by the equation
\beq\label{example}
w^2+Q_m(E)w-E^{2m+1}=0
\eeq
where $Q_m$ is a polynomial of degree $m$.
\end{rem}

\medskip
\noindent{\it Proof of Lemma \ref{lem-curve}.} By definition, the polynomial $R$ is of degree $k+1$ in $w$ and of degree $n$ in $E$. The matrix $L(w)$ has $(k+1)$ non-zero diagonals below the main diagonal and $(k+1)$ diagonals in the upper right corner. The entries of the latter are multiplied by $w$. Hence, $R(0,E)=E^n$, and the summation in (\ref{R}) is going over $i>0$. The coefficient $r_{1,\,0}$ is the product of the entries of the first above the main diagonal and the coefficient at $w$ in the lower left corner of the matrix. That proves (\ref{r}).

It remains only to prove that the summation in (\ref{R}) goes over the pairs of indices $(i,j)$ such that $ni+j(k+1)<n(k+1)$. For further use we present the proof of the last statement which goes along the lines used in the theory of commuting differential operators (see details in \cite{kr-nov}).
\begin{lem}\label{lem-psi} Let $L$ be in $\D_{k+1,\,n}$. If $n$ and $(k+1)$ are coprime, then there exists a unique formal series
\beq\label{e}
E(z)=z^{-(k+1)}\left(1+\sum_{s=1}^\infty e_s z^s\right)
\eeq
such that equation $L\psi=E\psi$ has the unique formal solution
\beq\label{psiinf}
\psi_i(z)=z^{i}\left(1+\sum_{s=1}^\infty \xi_s^+(i)z^{s}\right)
\eeq
with periodic coefficients
\beq\label{xiper}
\xi_s^+(i)=\xi_s^+(i+n)
\eeq
and normalized by the condition $\xi_s^+(0)=0$.

\end{lem}
\noindent {\it Proof.} The substitution of (\ref{psiinf}) and (\ref{e}) into the equation $L\psi=E\psi$ gives a system of difference equations for unknown coefficients $e_s$ and $\xi_s$ of the series.
The first of them is the equation
\beq\label{xi1}
e_1+\xi_1^+(i)-\xi_1^+(i-k-1)=a_i^{k}\,.
\eeq
The periodicity constraint (\ref{xiper}) uniquely defines $e_1=n^{-1}\sum_{i=1}^na_i^{k}$ and reduces difference
equation (\ref{xi1}) of order $k+1$ to the difference equation of order 1:
\beq\label{xi11}
me_s+\xi_1^+(i)-\xi_1^+(i-1)=\sum_{j=0}^{m-1} a_{i-j(k+1)}^{k}\,,
\eeq
where $m$ is the integer $1\leq m<n$ such that $m(k+1)=1 \, ({\rm mod}\ n)$. Equation (\ref{xi11}) and the initial condition $\xi_1^+(0)=0$ uniquely defines $\xi_1^+(i)$.

For arbitrary $s$ the defining equation for $e_s$ and $\xi_s^+$ has the form:
\beq\label{xi}
e_s+\xi_s^+(i)-\xi_s^+(i-k-1)=Q_s(e_1,\ldots, e_{s-1};\xi_1,\ldots,\xi_{s-1})\,.
\eeq
where $Q_s$ is an explicit function linear in $e_{s'}, \xi_{s'},\ s'<s,$ and polynomial in $a_i^j$.
The same arguments as above show that it has the unique periodic solution. The lemma is proved.

By definition the series $E(z)$ is formal. It turns out that it is a convergent series and is the Puiseuz series solution of equation (\ref{R}) at infinity. More precisely, the arguments identical to that in (\cite{kr-com}, see details in \cite{kr-nov}) prove that:
\begin{lem} Let $E(z)$ be the formal series defined in Lemma \ref{lem-psi}, then the characteristic polynomial
$R(w,E)$ is equal to
\beq\label{RE}
R(w,E)=\prod_{j=1}^n (E-E(z_j)) , \ \ z_j^{-n}=w
\eeq
\end{lem}
The right hand side of (\ref{RE}) is a symmetric function of the variables $z_j$. Hence a priori it is a polynomial in the variable $E$ and a (formal) Laurent series in $w^{-1}$. Equation (\ref{RE}) states that it is a polynomial in $w$. It also implies also that the degree in $z$ of all terms in $R(w,E)$ except $w^{k+1}$ and $E^n$ is strictly less then $n(k+1)$. Lemma \ref{lem-curve} is proved.

Similarly we describe the Bloch solution near the second marked point.

\begin{lem}\label{lem-psi-zero} If $L$ is a "strictly" lower triangular operator (not necessarily periodic), then the equation $L\psi=E\psi$ has a unique formal solution of the form
\beq\label{psizero}
\psi_i(E)=e^{x_i}E^{i}\left(1+\sum_{s=1}^\infty \xi_s^-(i)E^{s}\right), \ e^{x_i-x_{i-1}}:=a_i^1.
\eeq
and normalized by the condition $\xi_s^-(0)=0$.
\end{lem}
{\it Proof.} The substitution of (\ref{psizero}) into (\ref{eigenL}) gives a system of equations for unknown coefficients $\xi_s^-$. They are nonhomogeneous first order difference equations
\beq\label{chi}
\xi_s^-(i)-\xi_s^-(i-1)=q_s(\xi_1^-,\ldots,\xi^-_{s-1})\,,
\eeq
which together with the initial conditions recurrently define $\xi_s^-(i)$ for all $i$.

The uniqueness of the series solution (\ref{psizero}) implies
\begin{cor} If $L\in \D_{k+1,\,n}$, then the formal series (\ref{psizero}) is the Bloch solution, i.e. it satisfies (\ref{monod}) with
\beq\label{wzero}
w(E)=\psi_n(E)=E^n\sum_{s=0}^\infty w_s E^s.
\eeq
\end{cor}
The coordinates of the Bloch solution $\psi_i=\psi_i(p),\ p=(w,E)\in \G$, normalized by the condition $\psi_0$
are rational functions of $w,E$. By standard arguments the poles of $\psi_i$ do not depend on $i$. In order to find the degree of the pole divisor $D$ it is enough to consider the matrix $F(E)$ with entries $F^{ij}=\psi_i (p_j(E)),\ 0\leq i\leq k$, where $p_j(E)=(w_j,E),\ 1\leq j\leq k+1,$ are the preimages of $E$ on $\G$ under the projection map $p\in\G\to E(p)\in \mathbb C$.
The matrix $F$ depends on the ordering of the preimages but the function $f(E)=\det^2 F(E)$ does not, i.e. $f$ is a rational function of $E$. It has double poles at the projections of the poles of $\psi$ (if they are distinct) and the  pole of order $2k$ at $E=0$. If $\G$ is smooth then $f$ has zeros at the finite branching points of the cover $\G\to \bbC$. The multiplicity of the zero of $f$ equals to the multiplicity of the branching. The infinity point $p_+$ is the branching point of the cover $\G\to \bbC$ of multiplicity $k$. From (\ref{psiinf}) it follows that at $E=\infty$ the function $f$ has zero of order $k$. The degrees of the pole and the zero divisors of a rational function coincide. Hence, $2\deg D+2k=\nu$, where $\nu$ is the total multiplicity of the branch points. The cover $\G\to \bbC$ is of degree $k+1$. Hence, by the Hurwitz formula the genus of $\G$ (if it is smooth) is equal to $2g=\nu-2k$. Therefore, $\deg D=g$ and we have proved that the Bloch solutions on the smooth spectral curves coincide with the discrete Baker-Akhiezer function.

\begin{rem} The coefficients $r_{ij}$ of the characteristic equation (\ref{R}) are polynomial functions of the coefficients of the operators $L\in \D_{k+1,\,n}$. The direct proof that they are independent is highly nontrivial and the author is not aware of any universal approach to this problem except the combined use of the direct and the inverse spectral transforms. In order to prove their independence it is enough to show that the image of the correspondence $\D_{k+1,\,n}\to \SP$ contains at least one smooth spectral curve. The latter is established by the construction of the inverse map (\ref{correspondence}).
\end{rem}
Summarizing the results presented above we obtain the following statement.

\begin{theo} If $n$ and $k+1$ are coprime, then the map (\ref{correspondence}) is a one-to-one correspondence between an open subset of the Jacobian bundle over the family of curves given by (\ref{R}) and an open subset in the space $\D_{k+1,\,n}$ of $n$-periodic lower triangular difference operators of order $k+1$.
\end{theo}

\section{Superperiodic difference operators}

In this section we consider the spectral theory of superperiodic operators.
The spectral curve $\G_{spec}$ of an operator $L\in \E_{k+1,\,n}$ is never smooth. The definition of
$\E_{k+1,\,n}$ implies that the point $p_{\,0}=(w=(-1)^{n+k}, E=-1)\in \G_{spec}$ is a multiple point of order $(k+1)$ (at which all sheets of the cover $E:\G_{spec}\to \mathbb C$ intersect transversally). The latter requires vanishing of the coefficients of the Tailor expansion of $R$ at $p_0$ of degree less than $k$. That is
$\frac{(k+1)(k+2)}2$ linear equations on the coefficients of $R$. There is one relation between these equations because  if $w=(-1)^{n+k}$ is a root of the equation $R(w,-1)=0$ of multiplicity {\it at least} $k$, then it is of multiplicity $(k+1)$. Hence, the space $\SP_{spec}$ of the spectral curves of operators
$L\in \E_{k+1,\,n}$ is of dimension $\frac{k(n+1)}2-\frac{(k+1)(k+2)}2+1=\frac{k(n-k-2)}2$

Let $\Sigma_{k+1,\,n}$ be a space of algebraic curves $\G$ that are the partial normalization $\pi:\G\to \G_{spec}$ resolving the multiple point $p_0$, i.e. it is one-to-one except at $(k+1)$ smooth point $p^{\,j}\in \G,\ j=1,\ldots, k+1$, that are the preimages of $p_0$, i.e., $\pi(p^{\,j})=p_{\,0}$.
In other words the smooth curve $\G$ in $\Sigma_{k+1,\,n}$ is characterized by the following properties:

(i) {\it on $\G$ there exists a meromorphic function $w$ having one pole at $p_+$ of order $n$ and having zero of order $n$ at another point} $p_-$;

(ii) {\it on $\G$ there exists a meromorphic function $E$ having one pole at $p_+$ of order $k+1$ and such that
the zeros $p^j$ of the function $E+1$ are distinct and at $p^j$ the function $w-(-1)^{k+n}$ also vanishes}:
\beq\label{pset}
E(p^{\,j})=-1,\ \ w(p^{\,j})=(-1)^{k+n},\ \ j=1,\ldots,k+1\,.
\eeq

\noindent
Under the normalization $\pi$ the genus of $\G_{spec}$ decreases by $\frac{k(k+1)}2$. Hence, for an open set in $\Sigma_{k+1,\,n}$, where $p_{\,0}$ is the only singularity of $\G_{spec}$, the corresponding curve $\G$ is smooth of genus $g=\frac{k(n-1)}2-\frac{k(k+1)}2=\frac{k(n-k-2)}2$.

Consider the preimage $\pi^*(\psi)$ (for which we will keep that same notation $\psi$) of the Bloch function under the normalization map. It has the same expansions (\ref{psiinf},\ref{psizero}) at the marked points $p_{\pm}\in \G$. The function $E(p), p\in \G$ defines a cover $\G\to \mathbb C$ which has $(k+1)$ sheets over the point $E=-1$. The same counting as above gives that $\psi$ has $g$ poles on $\G$. Hence, we get the following statement:
\begin{theo} For an open set of operators $L\in \E_{k+1,\,n}$ the Bloch solutions $\psi$ of equation $L\psi=E\psi$ are parameterized by points $p$ of a smooth algebraic curve $\G\in \Sigma_{k+1,\,n}$. Moreover, the function $\psi_i(p)$ is the Baker-Akhiezer function.
\end{theo}

Notice, that the evaluation of formula (\ref{psitheta})  at the points $p^{\,j}\in \G$ gives explicit theta-functional expression of the basis $\Psi^{(j)}:=\psi(p^j)$ of solution of the equation $(L+1)\Psi=0$.

\begin{cor} The inverse and the direct spectral transforms establish a one-to-one correspondence
\beq\label{corrab}
\P\ \ \Leftrightarrow\ \  \E_{k+1,\,n}
\eeq
between an open set of the Jacobian bundle $\P$ over the subspace of smooth curves in $\Sigma_{k+1,\,n}$ and an open set of superperiodic operators.
\end{cor}
Now we are ready to prove Theorem \ref{main}.

{\it Proof.} Consider the function
\beq\label{G}
G=\frac{w-(-1)^{n+k}}{E+1}+(-1)^{n+k}
\eeq
on $\G\in \Sigma_{k+1,\,n}$. At the zeros $p^{\,j}$ of the denominator the numerator vanishes as well. Hence, $G(p)$ is holomorphic on $\G\setminus p_+$. At the marked point $p_+$ the function $G$ has pole of order $n-k-1$. It vanishes at $p_-$. By Lemma \ref{comring} there exists an operator $\L=L_G\in \D_{n-k-1,n}$ such that $\L\psi=G\psi$, where $\psi$ is the Bloch solutions of the equation $L\psi=E\psi$.
From the definition of $G$ it follows that operator equation
\beq\label{opereq}
T^{-n}-(-1)^{k+n}=(\L-(-1)^{k+n})(L+1)=(L+1)(\L-(-1)^{k+n})
\eeq
holds. The part $(i)$ of the theorem is proved.

The function $w-(-1)^{k+n}$ has $n$ zeros on $\G$. $(k+1)$ of them are the points $p^{\,j}$. Let $q^{\,m},\ m=1,\ldots, n-k-1$ be the set of the remaining zeros. At these points $G(q^m)=(-1)^{k+n}$. Hence, $\L\in \E_{n-k-1,n}$.

In order to prove the last statement of the theorem it remains only to prove equation (\ref{duality}). The latter is a direct corollary of the following statement:
\begin{lem} Let $\widehat \Psi_{(i)}$ be $(n\times(k+1))$ matrix with columns
\beq\label{widehatpsi}
\left(\psi_{i+1}(p^{\,j}),\ldots,\psi_{i+n}(p^{\,j})\right)
\eeq
and let $\widehat \Psi_{(i)}^+$ be $((n-k-1)\times n)$ matrices with rows
\beq\label{wideharpsi1}
\left(\psi_{i+1}^+(q^{\,m}),\ldots,\psi_{i+n}^+(q^{\,m})\right)
\eeq
where $\psi$ and $\psi^+$ are the Baker-Akhiezer function and its dual on $\G\in \Sigma_{k+1,\,n}$ and $\{p^j,q^m\}$ are zeros of the function $w-(-1)^{k+n}$. Then the orthogonality relation
\beq\label{ortog}
\widehat \Psi_{(i)}^+ \ \widehat \Psi_{(i)}=0.
\eeq
holds for all $i$.
\end{lem}
{\it Proof.} The relation (\ref{ortog}) is a corollary of the orthogonality of eigenvectors and covectors of
a linear operator corresponding to different eigenvalues.
\beq\label{ortog1}
E(q^{\,m})<\psi^+(q^m) \psi(p^{\,j}>=<\psi^+(q^m) L\psi(p^{\,j})>=E(p^j)<\psi^+(q^m)\psi(p^{\,j})>
\eeq
where $<\cdot,\cdot>$ stands for the pairing between vectors and covectors.
By definition $E(p^{\,j})=-1$. The function $E$ is of degree $k+1$. Hence, it does not equal to $-1$ at any other point of $\G$. That implies $E(p^{\,j})\neq E(q^m)$. Then, from (\ref{ortog1}) we get the equation
\beq\label{or1}
<\psi^+(q^m) \psi(p^{\,j})>=\sum_{s=1}^n \psi_{i+s}^+(q^m)\psi_{i+s}(p^{\,j})=0\,,
\eeq
which is equivalent to (\ref{ortog}).

For $L\in \E_{k+1,n}$ equation (\ref{ortog}) coincides with the duality relation (\ref{matrixdual}) between solutions of the corresponding equations for the operators $L=L_E$ and $\varepsilon L_G^\s \varepsilon^{-1}$.
The part $(ii)$ of Theorem \ref{main} is proved.

\bigskip
\noindent {\it Proof of Theorem \ref{main1}.} The direct spectral transform (\ref{corr}) is constructed along the same lines as in the case of the spectral theory of periodic operators. Let us briefly recall some necessary steps
of the construction. Let $\D_{k+1}$ be the affine space of monic strictly triangular operator of order $(k+1)$
(not necessary periodic). Consider a pair $L\in \D_{k+1}$ and $K\in \D_{n-k-1}$ of commuting operators, $[L,K]=0$.

The restriction of operator $K$ onto $(k+1)$-dimensional space of solution of the equation $(L-E)\psi=0$ defines a finite-dimensional linear operator $K(E)$.

\begin{lem}[\cite{kr-diff}]\label{com-curve} If $n$ and $k+1$ are coprime, then the characteristic equation of the operator $K(E)$ has the form:

\beq\label{Rnew}
\R(\kappa,E)=\kappa^{k+1}-E^{n-k-1}+\sum_{(i,j)\in \,I}\,\rho_{ij}\, \kappa^iE^j=0\,,
\eeq
where the summation is taken over the set of pairs of non-negative integers $(i,j)$ such that $0<i(n-k-1)+j(k+1)< (n-k-1)(k+1).$
For commuting pairs in general position (when the spectral curve $\G$ defined by (\ref{Rnew}) is smooth) the common eigenfunction of the commuting operators, $L\psi=E\psi, \ K\psi=\kappa\psi,\ (\kappa,E)\in \G$ is the discrete Baker-Akhiezer function.
\end{lem}
\noindent
The number of the coefficients in (\ref{Rnew}) is $\frac{(k+2)(n-k)}2-2$. For generic values of these coefficients the corresponding curve $\G$ is smooth and has genus $g=\frac{k(n-k-2)}2$\,. Two marked points on $\G$ are the infinity point $p_+$ where $E$ and $\kappa$ have poles of orders $(k+1)$ and $(n-k-1)$, respectively, and the point $p_-$ where two functions vanish. If $\G$ is smooth then the ring $\A_+$ of meromorphic functions on $\G$ having pole only at $p_+$ and vanishing at $p_-$ is generated by $E$ and $\kappa$. If $2k+2<n$, then the function $E$ can be characterized as the generator with the lowest order of pole at $p_+$. By this definition it is unique up to the multiplication by a constant $E'=c^{k+1}E$. The second generator $\kappa$ of the ring is unique up to the transformation
\beq\label{automor}
\kappa'=c^{n-k}\kappa-\sum_{s=1}^{[\frac{n}{k+1}]-1}\, \a_s \, E^s\,,
\eeq
which preserves the order of pole at $p_+$ and the form of equation (\ref{Rnew}).

As it was mentioned above the generic pairs of commuting operators have quasi-periodic coefficients. By assumption of the theorem the operator $L$ is $n$-periodic, i.e. $[L, T^{-n}]=0$. That implies that the ring $\A_+$ contains a function $w$ which has the pole of order $n$ and $p_+$ and zero of order $n$ at $p_-$. Because $A_+$ is generated by $E$ and $\kappa$, this function can be represented in the form
\beq\label{wek}
w=\kappa E - e\kappa + Q(E)
\eeq
where $e$ is a constant and $Q$ is a polynomial of degree $d$ such that $d(k+1)<n$ (recall that by assumption $2(k+1)<n$) and vanishing at $E=0$.

Equation (\ref{wek}) implies that the function $w$ has the same value $\mu:=Q(e)$ at zeros $p^{\,j}$ of the function $E-e$. In general position these zeros are simple. Therefore, the evaluation of the Baker-Akhiezer at these points gives a basis of solutions of the equation $(L-e)\Psi=0$ having the same monodromy multiplicator $\mu$.
\begin{rem} The assumption that $\{p^{\,j}\}$ is a set of $(k+1)$ distinct points is the precise meaning of "general position" in the formulation of the theorem.
\end{rem}
The scaling of $E$ and $\kappa$ by a constant $c$ corresponds to the scaling transformation (\ref{scaling}) of operators. By proper choice of the scaling constant $c$ the constant $e$ in (\ref{wek}) can be transformed to $e=-1$. Hence, the operator $L$  up to the scaling transform is superperiodic.  The first statement of the theorem is proved.

For the proof of the second statement it  is enough to note, that the polynomial $Q(E)-Q(e)$ is always divisible by $(E-e)$, i.e. there is a unique polynomial $P(E)$ without the free term, $P(0)=0$, such that the equation $Q(E)-\mu=(E-e) (P(E)+e^{-1}\mu)$ holds. Then  from (\ref{wek}) we get the equation
\beq\label{wek1}
w(p)-\mu=(E(p)-e)(\lambda(p) +e^{-1}{\mu}), \ \ \lambda(p) =\kappa+P(E)
\eeq
The latter implies that the operator $\L=K+P(L)$ is superperiodic up to the scaling transform. The theorem is proved.

\medskip
\begin{rem} Throughout the paper we have considered mostly operators with complex coefficients, but all the constructions and results admit real reduction. If the coefficients of the operator $L$ are real, then the coefficients of the characteristic polynomials are real. Therefore, the complex conjugation defines a antiholomorphic involution $\tau$ of the corresponding spectral curves. Equation (\ref{tau}) implies that the divisor $D$ is invariant under the involution $\tau(D)=D$. Hence, the direct and the inverse spectral transform establish one-to-one correspondence between the space of operators with real coefficients and a bundle over the space of {\it real} spectral curves whose fiber is real locus of the corresponding Jacobian.
Note, that the locus of $M$-curves and divisors $D$ as in Lemma \ref{mcurves} corresponds
to operators $L$ with sign definite coefficients $a_j^1$.
\end{rem}

\end{document}